Low Progress Math in a High Performing System:

The Role of Math Anxiety in Singapore's Elementary Learners

# Introduction

Mathematics is indivisible from daily life, not only because it accounts for a large proportion of academic achievement in schools, but also because it involves in many areas such as future professional development, investment and financial management (Bynner, 2004; Lusardi & Mitchell, 2011; Mammarella et al., 2015). As a country with high-performing achievement in math, Singaporean students were found to show a higher level of math anxiety (Foley et al., 2017), which is defined as feelings of tension or worry that interfere with performance in math in a variety of daily activities and school settings (Richardson & Suinn, 1972), even in high performing education systems. In Singapore, while findings from the Programme for International Student Assessment (PISA) show that the city state has been consistently ranked as one of the world's top education systems particularly in math, there is also a long tail in achievement distributions (OECD, 2013). Specifically, there is still a 15-20 percentile band of low progress mathematics learners who continue to struggle despite multi-pronged behavioural intervention approaches in school (Ng et al., 2007). Some of them are also at risk of developmental dyscalculia (DD), a brain-based learning disability of mathematics which cannot be ascribed to low intelligence or/and inadequate education (Rubinsten & Henik, 2009). This poor initial math ability may be compounded by math anxiety, creating a vicious cycle (Rubinsten & Tannock, 2010). Yet, extant literature on math anxiety have been dominated by adults' research, ignoring the emergence of math anxiety in children, especially for those with poor math performances, such as DD at-risk.

To fill this research gap, we recruited low progress (LP) learners in math from Singapore who has been screened by a national math test, specifically the Ministry of Education's (MOE) Early Numeracy Indicator (ENI) which classifies students into either ENI 0, 1, 2, or 3 with ENI 0 being the weakest and ENI 3 the strongest group, to investigate the relationship between math anxiety and math performance. Meanwhile, as a multinational state, we also examine ethnic differences (Chinese VS Malay VS Indian VS Others) in math anxiety, which remains understudied in the context of Asian education systems. Considering the debating conclusions in gender difference (Szczygiel, 2020; Van Mier, 2019), gender was also included in the demographic characteristics.

## Literature Review

### Math Anxiety: A Function of Demographic Characteristics

From research on typical developing children, the influence of math anxiety has been found to differ based on demographic characteristics, such as gender and ethnicity. While limited research has investigated the differences of demographic characteristics in math anxiety on typically developing math achievers in young children (Szczygiel, 2020; Van Mier, 2019), Szczygiel (2020) found that girls showed higher math anxiety than boys in testing and total dimension of Abbreviated Math Anxiety Scale for Elementary Children (mAMAS-E, Szczygieł, 2019). However, no gender differences in the learning dimension of math anxiety in second graders (Szczygiel, 2020) were found, while Van Mier (2019) found that only girls showed significant math anxiety which negatively moderated math performance. These inconsistent conclusions formed our impetus of investigating gender differences when studying the effect of math anxiety in our participants' profile. Concurrently, the literature remains scant on the effects of math anxiety for learners with math learning difficulties (MLD), who, as argued by Dowker (2005) and Krinzinger (2009) are more vulnerable towards anxiety debilitative experiences. Although there was one research which investigated math anxiety in children with/without MLD, gender differences were not compared in children with MLD, which may be due to the small sample size (20 boys versus 9 girls; Mutlu, 2019). Another important demographic characteristic is ethnicity, where its relation to math anxiety remains understudied (Hart & Ganley, 2019; Ahmed, 2018). From a recent study on adults, Hart and Ganley (2019) reported that there was no ethnic group difference among Hispanic, Black and White adults. However, in studies on adolescents, Hispanics were found to have higher levels of math anxiety than their White counterparts, but Blacks did not differ from Whites (Cheema & Sheridan, 2015), while Ahmed (2018) found that Hispanics and Black adolescents showed a consistently high math anxiety trajectory compared to their White counterparts through a six-year longitudinal observation. In the current limited research, only Hispanic, Blank and White groups were investigated, no prior research has investigated ethnic group differences in math anxiety in Asian contexts, such as Malay, Indian, Chinese and others, especially in young children.

**Working Memory: The Mediator in Math Anxiety and Math Performance**

Math anxiety has been found to be consistently related to poor math performance (Ashcraft, 2002; Dowker et al., 2016; Hembree, 1990; Namkung et al., 2019; Zhang et al., 2019). Longitudinal studies suggest that math anxiety is present and stable even from a young age and throughout adulthood (Cargnelutti et al., 2017; Krinzinger et al., 2009; Sorvo et al., 2019; Vukovic et al., 2013). Yet how math anxiety negatively affects mathematics performance is still debated. One of the most studied factors that could account for this negative relationship is working memory (Justicia‐Galiano et al., 2017). Working memory is a limited mental workspace that handles the storage of information and manipulates it simultaneously (Baddeley & Hitch, 1974;

Justicia-Galiano et al., 2017). Literature has shown that working memory ability is involved in mathematics-related tasks (see Raghubar, Barnes & Hecht, 2010 for a review). For example, in a simple one-digit addition task, young children have to remember all the addends and operator and then manipulate them simultaneously. According to the attentional control theory (Eysenck et al., 2007), the worrying thoughts and feelings related to mathematics-related tasks would reduce the storage capacity and hamper information processing within limited working memory workspace. In one study on children aged 8–12 years, multiple mediation analyses indicated that working memory, as a mediator, contributed to explaining the relationship between math anxiety and math performance (Justicia-Galiano et al., 2017). A further research on children aged 5-8 years found that a negative relation between math anxiety and math performance was only found in children with higher working memory. The possible reason could be that children with high working memory tend to employ strategies requiring intensive working memory resources (Ramirez et al, 2013). Findings from neuroimaging studies also provide plausible evidence of how working memory is disrupted in the negative association between math anxiety and math performance. In a neuroimaging study on children with math anxiety, Young, Wu and Menon (2012) found that children with higher math anxiety showed higher neural responses in the amygdala generating negative emotions while lower response in the dorsolateral prefrontal cortex subserving working memory. To deal with the aroused math anxiety, individuals had to distribute attention source from the central executive network that is engaged during goal-directed processing, thus leaving less workspace for working memory to process math tasks efficiently (Pletzer et al., 2015; Menon & Uddin, 2010). This may indicate that the aroused math anxiety disrupts the working memory required for processing math tasks.

**Dyscalculia: The Moderator in Math Anxiety and Working Memory**

From the summarized literature, it is known that math anxiety may disrupt working memory in typically developing learners. But will the disruption be exacerbated in low progress math learners and those at-risk of developmental dyscalculia (DD)? In reviewing literature related to math anxiety and working memory, firstly it has been reported that math anxiety is often associated with low math performance, such as mathematics learning difficulties or even developmental dyscalculia. In a large-scale study which examined 1757 primary- and secondary-school children, it was found that children with DD were twice as likely to have high math anxiety as children with typical mathematics performance (Devine et al., 2018). Second, working memory disruption, especially in the visuospatial working memory, is also found in children with DD. Behavioural research found that the visuospatial working memory was disrupted in children with DD, compared with matched typically developing children on general intelligence, reading and other cognitive measures (Ashkenazi, 2013) or children with reading difficulties (Lee, Howard & Saez; 2006), see a review from (Menon, 2016). In a neuroimaging study, children with DD were found to show weaker neural activation during a spatial working memory task in the right intraparietal sulcus (IPS) in a functional magnetic resonance imaging (fMRI) scanner and impaired working memory proficiency outside the scanner, compared to typically developing children. Additionally, significant correlations of the right IPS activity with performance on the working memory task were also found (Rotzer, 2009). Thirdly, the disruptive role of math anxiety in working memory may indicate a negative relationship between both. Evidence from Ashcraft and Kirk (2001) showed that adults with high math anxiety performed worse in working memory tasks than adults with low or medium math anxiety. Evidence from children showed similar results. Mammarella et al. (2015) examined the relationship between math anxiety and working memory in sixth- to eight-grader children and found that children with math anxiety were particularly impaired in the working memory task, compared to children without math anxiety. Taken together, considering the complex relationship between math anxiety and working memory in typically developing, low progressing or children with DD, we posit a plausible hypothesis that DD may moderate the relationship between math anxiety and working memory such that children with dyscalculia may show a stronger negative relationship between math anxiety and working memory. To this end, this paper reports our investigations on the following strands of inquiry:

RQ1) Does math anxiety differ significantly between demographic characteristics specifically gender, ethnicity?

RQ2) Does math anxiety explain significant variance in predicting math performance after controlling for related covariables such as demographic characteristics and working memory?

RQ3) Does math anxiety level differ in relation to working memory capacity?

RQ4) Does dyscalculia-at-risk status moderate the interaction effect between math anxiety and working memory capacity?

## Method

**Participants**

Participants ($N = 151$) were first graders from 17 Singapore elementary schools identified as low math achievers through the ENI and who were participants of a larger mathematics study (Jamaludin & Hung, 2019). The participants consisted of 93 males and 58 females with 49.67% of the participants of Malay ethnicity, 39.07% were Chinese, and 6.62% were Indian. The participants' age ranged from 74.78 to 100.81 months with mean of 82.34 months ($SD = 4.27$). No learning or psychological diagnosis were reported during the study registration among the participants.

Upon entry into Grade 1, participants had completed a national math screening test, the Early Numeracy Indicator (ENI), and were grouped into four different levels (i.e. ENI0, 1, 2 and

3) based on their scores. Typically, children identified as having ENI score of 2 and below will need to attend a pull-out math program for two years called the Learning Support for Mathematics (LSM), although selected ENI 3 children might also be assigned to the LSM program based on teachers' recommendation. All children participants in the present study were from the LSM program and consisted of 26 children identified as ENI0 (lowest possible level), 61 as ENI1, 35 as ENI2, and 29 as ENI3.

**Measures**

*Math Anxiety*

Math anxiety was measured using Scale of Early Math Anxiety (SEMA; Wu et al., 2012). SEMA is a 20-item questionnaire which requests children to report their level of nervousness on a five-point Likert's scale when presented with math related questions or situations. The questionnaire is divided into two sections with the first 10 items assessing anxiety related to Numerical Processing (e.g. "Is this right? 9+7=18") and the last 10 items related to contextual situations related to math (e.g. "You are about to take a math test"). The total score ranged from 0 to 80 with higher score indicating higher math anxiety.

*Working Memory*

Working memory in this study was measured using Digit Span and Picture Span subtests of the *Weschler Intelligence Scale for Children-Fifth Edition* (WISC-V; Wechsler et al., 2014). The Digit Span consists of three tasks which requires the child to repeat numbers read aloud by the presenter exactly as spoken (Forward), in backward order (Backward), and in specific sequence (Sequential). In Picture Span, the child is required to recall the order of a set of pictures presented. The raw scores of the subtests were summed and age-scaled to create the Working Memory Index (WMI).

*Verbal Ability*

Considering high correlation between math and verbal ability (Henry, Baltes & Nistor, 2014), verbal ability was measured as a covariable. It was measured using Similarities and Vocabulary subtests of the *Weschler Intelligence Scale for Children-Fifth Edition* (WISC-V; Wechsler et al., 2014). The Similarities subtest requires the child to describe how two words are similar while the Vocabulary subtest requires the child to define words. The raw scores of the subtests were summed and age-scaled to create the Verbal Comprehension Index (VCI).

*Math Performance*

Math performance was measured using the Math Composite score of the *Wechsler Individual Achievement Test-Third Edition* (WIAT-III; Psychological Corporation, 2009) which consists of Numeracy and Math Problem Solving subtests. The Numeracy subtest is a paper-and-pencil test that measures the ability in basic math skills and basic math operations such as integers, geometry, and algebra. Math problem solving subtest requires the child to provide verbal or pointing responses on questions measuring the child ability in their understanding of basic math concepts, everyday applications, geometry, and algebra. The raw scores of the two tests were summed and age-scaled to create the Math Composite score.

*Dyscalculia Screener*

As part of the larger study, participants were screened using the Dyscalculia Screener (Butterworth, 2003) which consisted of four tasks: Simple Reaction Task, Dot Enumeration, Numerical Stroop, and Addition test. The tasks are computerized and administered using a laptop. The child follows the instruction provided on screen using a headphone and responds to the stimuli by pressing buttons on the right or the left side of the keyboard. The raw scores of the Dot Enumeration, Numerical Stroop, and Addition test were standardized and were ranked on a range of 1 to 9. Based on the scores, a report of possible diagnosis is generated to indicate whether the child is at-risk of dyscalculia or not based on UK normed data. Importantly we only flag out dyscalculia *at-risk* status for our respective lines of inquiry – no formal DD diagnoses were made.

**Procedure**

All participants were part of a larger math intervention study cohort, iMAGINE (Jamaludin & Hung, 2019; NIE Office of Education Research, 2022), which aimed to understand the neurodevelopmental basis of math for children identified as having math learning difficulties. All the measures used in the current study were administered prior to any intervention done. SEMA and the Dyscalculia Screener were administered prior to the administration of WISC-V and WIAT-III. All the measures were administered one-to-one by trained research assistants at the children's school. Questions were read to and explained by the research assistants for all participants, so that those who had reading difficulties are not disadvantaged in any way.

*Data Analysis Procedure*

Missing data analysis was conducted on SPSS version 23 to test the Missing Completely at Random hypothesis (MCAR; Little, 1988). Result of Little's MCAR test showed that the hypothesis of missing completely at random could not be rejected, $\chi^2(50) = 53.98$, $p = .325$ indicating that the missing data were missing completely at random. Missing data were handled using Expectation-Maximization imputation using Amelia-II package version 1.7.6 (Honaker et al., 2019). Multivariate outlier analysis conducted using Mahalanobis' distance found one multivariate outlier. Analysis with and without the outlier did not change the findings and thus the outlier was retained in the analysis. The descriptive of the data after imputation is as shown in Table 1. All the hypotheses testing were conducted using R version 4.0.2 (R Core Team, 2020).

**Findings**

**Math Anxiety and Demographic Characteristics**

RQ1) Does math anxiety differ significantly between demographic characteristics specifically gender, ethnicity?

Type II ANOVA was used to investigate whether math anxiety measured using SEMA differed by gender, ethnicity, ENI scores, and suspected dyscalculia status (dyscalculia at-risk). Type II ANOVA was used to control for family-wise alpha error and was also found to be robust for unbalanced group sample sizes (Langsrud, 2003). Results showed that math anxiety did not vary significantly among the demographic characteristics, see table 2. Post-hoc tests using Tukey HSD found that the mean math anxiety score for children identified as ENI3 (M = 16.1, SD = 14.7) significantly different from children identified as ENI0 (M = 29.8, SD = 17.6) indicating that children identified as better math performer based on ENI test (ENI 3) showed lower math anxiety than the poorest performer on the ENI test (ENI 0) after controlling for family-wise alpha errors.

**Math Anxiety and Math Performance**

RQ2) Does math anxiety explain significant variance in predicting math performance after controlling for related covariables such as demographic characteristics and working memory?

To test the hypothesis of the effect of math anxiety on math performance, a hierarchical linear regression with the covariates consisting of gender, ethnicity, WISC-V's Verbal Comprehension Index (VCI), ENI levels and dyscalculia at-risk status, were entered in Block 1 followed by math anxiety and WISC-V's Working Memory Index (WM) in Block 2, see table 3. Result showed that the addition of Math Anxiety significantly improved the model $\Delta R^2 = 0.09$, $\Delta F(1, 139) = 5.38$, $p = .022$. Block 2 explained 40.73% of the total variance, $_{adjusted}R^2 = 0.36$, $F(11,$

139) = 8.68, $p < .001$. Math anxiety significantly predicted Math Performance, $b = -0.011$, $t(139) = -2.32$, $p = 0.022$, $_{partial}\eta^2 = 0.04$ such that lower math anxiety predicted higher math performance.

Among the covariates, language ability measured using WISC-V's VCI significantly predicted math performance such that higher verbal comprehension index predicted higher math performance, $b = 0.22$, $t(139) = 4.19$, $p < 0.001$, $_{partial}\eta^2 = 0.23$. Suspected dyscalculia status significantly predicted performance such that those suspected of dyscalculia performed significantly poorer at math compared to those not suspected of dyscalculia, $b = -0.62$, $t(139) = -3.44$, $p = 0.001$, $_{partial}\eta^2 = 0.1$. Working memory also significantly predicted math performance such that higher working memory predicted higher math performance, $b = 0.2$, $t(139) = 3.73$, $p < .001$, $_{partial}\eta^2 = 0.04$.

**Math Anxiety and Working Memory**

RQ3) Does math anxiety level differ in relation to working memory capacity?

To test the hypothesis of the interaction effect between math anxiety and working memory on math performance, the interaction effect was added into the Block 2 model. Results showed that the addition of the interaction effect did not significantly improve the model, $\Delta R^2 = 0.001$, $\Delta F(1, 138) = 0.13$, $p = .720$. There was no significant interaction effect between math anxiety and working memory.

**Moderating Effect of Dyscalculia 'at-risk' Status**

RQ4) Does dyscalculia-at-risk status moderate the interaction effect between math anxiety and working memory capacity?

To test the hypothesis of the moderating effect of suspected dyscalculia status on the interaction effect between math anxiety and working memory on math performance, the interaction effects were added into the Block 3 model. Result showed that the addition of the interaction effects did not significantly improve the model, $\Delta R^2 = 0.021$, $\Delta F(3, 135) = 1.64$, $p = .180$; indicating that there was no moderating effect of dyscalculia at-risk status on the interaction effect between math anxiety and working memory. Closer look at the regression coefficient however, showed that there was a significant 3-way interaction effect between those suspected as having dyscalculia, math anxiety, and working memory compared to those not suspected as having dyscalculia, $b = 0.014$, $t(135) = 2.04$, $p = 0.044$.

Analysis of the 3-way interaction using estimated marginal means of linear trend (Lenth et al., 2020) showed that except for normative children with average working memory, $b = -0.12$, $t(135) = -2.07$, $p = 0.04$, and high working memory, $b = -0.22$, $t(135) = -2.78$, $p = 0.01$, the simple slopes were not significantly different than zero indicating that math performance did not vary significantly as math anxiety increases for different levels of working memory except for normative children with average and high working memory which showed significant decreased in math performance as math anxiety increased, see figure 1.

Closer look at the interaction plot of the 3-way interaction showed negative slope of math performance as math anxiety increased except for children suspected for dyscalculia, $b = 0.09$, $t(135) = 0.71$, $p = 0.48$, this slope however was not significantly different from zero. Pairwise comparison of suspected dyscalculia and normative children with high working memory controlling for all family-wise errors was significant, $b = -0.31$, $t(135) = -2.07$, $p = 0.04$ indicating the slopes were significantly different from one another. No other pairwise comparisons were significant.

**Discussion**

This study examined math anxiety in a specific sample of children identified as having math learning difficulties (MLD) by their teachers. Previous research suggests that math anxiety is present for very young children, could be affected by various demographic characteristics (Zhang et al., 2019), and that it affects children with higher working memory more than those with lower working memory (Namkung et al., 2019; Ramirez et al., 2013; Vukovic et al., 2013).

The current study found that math anxiety did not vary significantly among the demographic characteristics except for those with very low math ability (ENI 0) and age-typical math ability (ENI 3). Our findings further found similar negative relationship between math anxiety and math performance in the sample of LP children with MLD. Similar relationship between math anxiety and working memory was also found such that those with average and high working memory performed poorer in math at higher levels of math anxiety but this was only true for those not suspected of dyscalculia. This supported the depleted cognitive capacity theory such that the cognitive capacity for those with higher working memory could be overloaded with worry which in turn affects math performance (Beilock & Carr, 2005; Kellogg et al., 1999; Mattarella-Micke et al., 2011).

**Demographic Characteristics of Math Anxiety**

As the findings indicate, math anxiety did not vary significantly between analyzed demographic characteristics among children identified as MLD except for the those with the weakest math ability (ENI 0) who had higher mean anxiety levels compared to age-typical math ability students on the ENI test (ENI 3). The present study's sample present a unique situation where initially only children with ENI score 2 and below were expected to be part of the sample. Discussions with the grade teachers indicated that while the national guideline was to use ENI 2 and below as indicator of MLD, some children identified as ENI 3 were also pulled out from the classroom for the math intervention program based on teachers' recommendation. This could be based on various reasons that are typically observed in Singapore's elementary mathematics classrooms as shared by teachers such as, weaker language abilities which might hinder them from following the math instructions in class. Closer look at the verbal reasoning index among ENI 3 children showed that on average, they performed below the population mean ($M = 75.52$, $SD = 13.11$) suggesting that these children might not have MLD but faced language issues compared to typically developing children instead. This paints a complex picture on how children could be identified as MLD which could be implicated by difficulties in language abilities.

**Complex nature of Math Learning Difficulty**

Results from the current study supported the negative link between math anxiety and math performance such that higher math anxiety predicted poorer math performance even after controlling for the children' language ability, working memory, gender, ethnicity, and math ability measured by the ENI scores and suspected dyscalculia status which have been shown to affect math performance (Geary et al., 2009; Vukovic et al., 2013). This finding adds to a larger body of work which found the debilitating effect of math anxiety in young children and especially for children identified as having MLD (Harari et al., 2013; Sorvo et al., 2019; Vukovic et al., 2013).

Additionally, the current study further found that among the covariates, language ability, working memory, and suspected dyscalculia status predicted math performance. This is consistent with other literature on the role of language (McClelland et al., 2007; Smedt et al., 2010) and working memory (Bull et al., 2008; Peng et al., 2016) on math performance.

The current study also found that those dyscalculia-at-risk status predicted poorer WIAT-III math composite score compared to those not suspected of dyscalculia based on the Dyscalculia Screener (Butterworth, 2003). The screener identifies children performing very poorly on core numerical tasks namely dot enumeration, magnitude comparison, and math knowledge compared to similar age children to determine the possibility of developmental dyscalculia. These tasks which have been shown to underlie math abilities indicated that these children had very poor foundational math skills and are at risk of having developmental dyscalculia (Gray & Reeve, 2016;

Reeve et al., 2012; Vetter et al., 2011). Overall, these findings suggest a complex nature of MLD which could be affected by multiple variables such as language, working memory, core numerical skills, and emotional factors.

Similar to the Vukovic et al. (2013) and Ramirez et al. (2013) study, the current findings also found that math anxiety disproportionately affected those with higher working memory. This relationship however was only present for children not suspected of dyscalculia. This raised an important question: why higher math anxiety is disproportionately more debilitating for those with higher working memory capacity but without dyscalculia at-risk status?

One possible explanation could be based on the depleted cognitive resource theory and poor math strategy utilization (Ramirez et al., 2016, 2018; Vukovic et al., 2013). It is possible that poor core numeracy skills i.e. enumerating and magnitude comparison which underlie those suspected of dyscalculia preclude effective utilisation of math problem solving strategies. This can be seen by consistently lower math performance at every level of working memory capacity for those identified as dyscalculia at-risk. Higher math anxiety thus might not lead to further decrease in math performance for those already faced with poor numerical skills. This is consistent with Ramirez et al.'s (2016) study which showed that poorer performance on math could be due to avoidance of advanced math problem solving strategies and reliance of simpler strategies instead.

**Limitations and Future Direction**

We discuss limitations in the current study that presents opportunity for further research in this area. Firstly, while the dyscalculia screener was useful in identifying children with poor core numeracy skills, the present study only included children already identified as 'low progress' math learners. It is not known whether there could be moderating effects between children not identified as MLD and suspected dyscalculia status. This could be useful to better understand how math anxiety affect working memory non-MLD, MLD, and MLD with poor core numeracy skills (developmental dyscalculia). Secondly, the current study did not include language ability as a potential moderator. Future studies could investigate alternative models to account for the complex nature of MLD which consist of unique interaction between the affect, cognition, and language abilities.

The differential impact of math anxiety on math performance for children with average and high working memory found in this study for children not suspected of dyscalculia only highlighted the complex nature of math learning difficulties in young children. Children identified as MLD in the current sample showed a myriad of emotional and cognitive difficulties including numeracy skills, math anxiety, and cognitive skills. Multi-dimensionality person-centred analyses aimed at identifying the various sub-types of MLD should be conducted to better understand and serve these vulnerable children who might be easily overlooked in a paradoxically high performing system (Yeo, Tan, Jamaludin, 2022). The aim is to provide more targeted and effective math support system *early* in their education journey that can attenuate otherwise accumulative learning struggles with mathematics over time.

**Acknowledgement**

This research was funded by the Singapore National Research Foundation (NRF) under the Science of Learning Initiative (NRF2016-SOL002-003).

# Tables and Figures

**Table 1.**

*Descriptive Statistics*

| Statistic (N = 151) | Mean | St. Dev. | Min | Max |
|---|---|---|---|---|
| Continuous | | | | |
| WIAT-III Composite Score | 87.68 | 12.07 | 60 | 120 |
| WISC-V WMI | 86.12 | 16.16 | 9 | 125 |
| WISC-V VCI | 74.18 | 16.52 | 35 | 108 |
| Math Anxiety | 23.71 | 17.32 | 0 | 80 |
| Age | 82.33 | 4.27 | 74.78 | 100.81 |
| Frequency | | | | |
| Gender | Male : 93 | Female : 58 | | |
| Ethnicity | Malay : 75 | Chinese : 59 | Indian : 10 | Others : 7 |
| ENI | ENI0 : 26 | ENI1 : 61 | ENI2 : 35 | ENI3 : 29 |
| Dyscalculia at-risk | No : 104 | Yes : 47 | | |

**Table 2.**

*ANOVA results for the relationship between demographic characteristics and math anxiety*

| | Sum Sq | df | $F$ value |
|---|---|---|---|
| Gender | 4 | 1 | 1.27 |
| Ethnicity | 5 | 3 | 0.54 |
| ENI | 21 | 3 | 2.44 |
| Suspected Dyscalculia | 6 | 1 | 2.18 |
| Residuals | 407 | 142 | |

**Table 3**

*Regression Coefficients for Multiple Regression Analyses with WIAT-III Math Composite score as the dependent variable*

|  | Block 1 | Block 2 | Block 3 | Block 4 |
|---|---|---|---|---|
| (Intercept) | 6.609 (0.435) *** | 5.678 (0.559) *** | 5.424 (0.904) *** | 4.495 (1.097) *** |
| Male (ref. Female) | 0.079 (0.179) | 0.145 (0.172) | 0.151 (0.173) | 0.130 (0.172) |
| Chinese (ref. Malay) | -0.126 (0.190) | -0.108 (0.179) | -0.116 (0.181) | -0.047 (0.182) |
| Indian | 0.253 (0.355) | 0.114 (0.335) | 0.099 (0.339) | 0.066 (0.338) |
| Others | 0.152 (0.430) | 0.311 (0.411) | 0.306 (0.412) | 0.293 (0.410) |
| WISC-V VCI | 0.281 (0.054) *** | 0.221 (0.053) *** | 0.221 (0.053) *** | 0.223 (0.054) *** |
| ENI1 | 0.269 (0.252) | 0.126 (0.239) | 0.139 (0.242) | 0.165 (0.244) |
| ENI2 | 0.162 (0.290) | -0.086 (0.278) | -0.074 (0.281) | -0.076 (0.283) |
| ENI3 | 0.651 (0.303) * | 0.306 (0.295) | 0.323 (0.300) | 0.350 (0.302) |
| Dyscalculia (DYS) | -0.715 (0.191) *** | -0.623 (0.181) *** | -0.627 (0.182) *** | 2.143 (1.851) |
| Math Anxiety (MA) |  | -0.011 (0.005) * | -0.001 (0.028) | 0.040 (0.038) |
| WISC-V WMI |  | 0.203 (0.055) *** | 0.231 (0.095) * | 0.336 (0.113) ** |
| MAxWM |  |  | -0.001 (0.003) | -0.006 (0.004) |
| MAxDYS |  |  |  | -0.108 (0.058) |
| WMxDYS |  |  |  | -0.343 (0.208) |
| MAxWMxDYS |  |  |  | 0.014 (0.007) * |
| R^2 | 0.384 | 0.407 | 0.408 | 0.429 |
| Adj. R^2 | 0.340 | 0.360 | 0.356 | 0.365 |
| Num. obs. | 151 | 151 | 151 | 151 |

**Figure 1.**

*3-way Interaction effect between Suspected Dyscalculia (dyscalculia at-risk) and Math Anxiety at different WISC-V WMI levels against Math Performance*

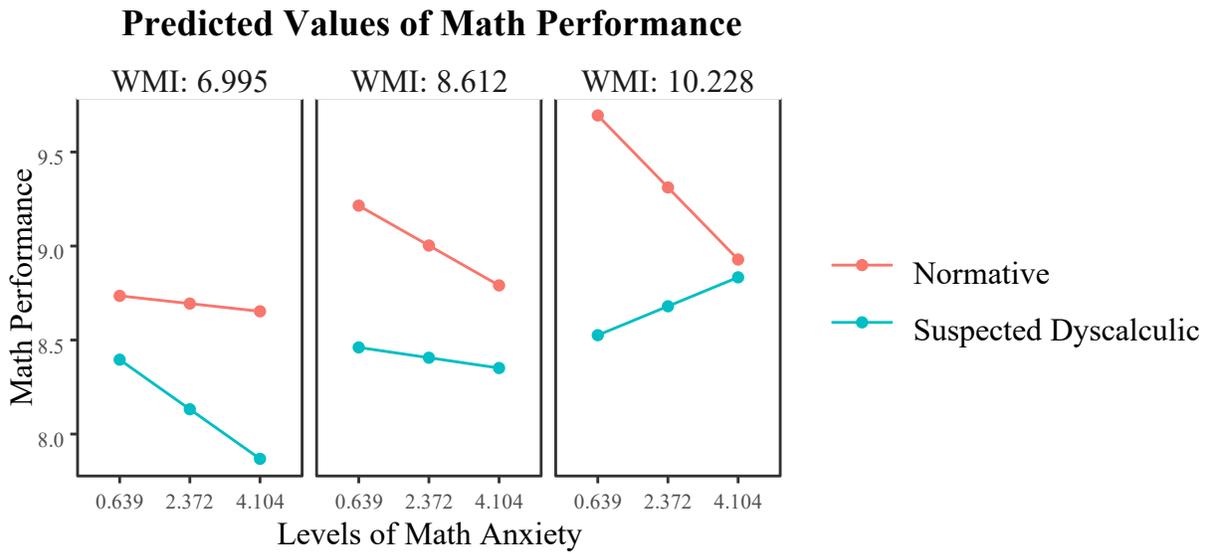

# References

Ahmed, W. (2018). Developmental trajectories of math anxiety during adolescence: Associations with STEM career choice. Journal of Adolescence, 67, 158-166.

Ashcraft, M. H., & Kirk, E. P. (2001). The relationships among working memory, math anxiety, and performance. Jounal of Experimental Psychology: General, 130, 224–237.

Ashkenazi, S., Rosenberg-Lee, M., Metcalfe, A. W., Swigart, A. G., & Menon, V. (2013). Visuo–spatial working memory is an important source of domain-general vulnerability in the development of arithmetic cognition. Neuropsychologia, 51(11), 2305-2317.

Ashcraft, M. H. (2002). Math anxiety: Personal, educational, and cognitive consequences. *Current Directions in Psychological Science*, *11*(5), 181–185. https://doi.org/10.1111/1467-8721.00196

Ashcraft, M. H., & Kirk, E. P. (2001). The relationships among working memory, math anxiety, and performance. *Journal of Experimental Psychology: General*, *130*(2), 224–237. https://doi.org/10.1037/0096-3445.130.2.224

Baddeley, A. D., & Hitch, G. (1974). Working memory. The Psychology of Learning and Motivation, 8, 47–89.

Bynner, J. (2004). Literacy, numeracy and employability: evidence form the British birth cohort studies. Literacy and numeracy studies, 13(1), 31-48.

Beilock, S. L., & Carr, T. H. (2005). When high-powered people fail: Working memory and "choking under pressure" in math. *Psychological Science*, *16*(2), 101–105. https://doi.org/10.1111/j.0956-7976.2005.00789.x

Bull, R., Espy, K. A., & Wiebe, S. A. (2008). Short-Term Memory, Working Memory, and Executive Functioning in Preschoolers: Longitudinal Predictors of Mathematical

Achievement at Age 7 Years. *Developmental Neuropsychology*, *33*(3), 205–228. https://doi.org/10.1080/87565640801982312

Butterworth, B. (2003). *Dyscalculia Screener* [Data set]. American Psychological Association. https://doi.org/10.1037/t05204-000

Carey, E., Hill, F., Devine, A., & Szűcs, D. (2016). The Chicken or the Egg? The Direction of the Relationship Between Mathematics Anxiety and Mathematics Performance. *Frontiers in Psychology*, *6*. https://doi.org/10.3389/fpsyg.2015.01987

Cargnelutti, E., Tomasetto, C., & Passolunghi, M. C. (2017). How is anxiety related to math performance in young students? A longitudinal study of Grade 2 to Grade 3 children. *Cognition and Emotion*, *31*(4), 755–764. https://doi.org/10.1080/02699931.2016.1147421

Chang, H., & Beilock, S. L. (2016). The math anxiety-math performance link and its relation to individual and environmental factors: A review of current behavioral and psychophysiological research. *Current Opinion in Behavioral Sciences*, *10*, 33–38. https://doi.org/10.1016/j.cobeha.2016.04.011

Cheema, J. R., & Sheridan, K. (2015). Time spent on homework, mathematics anxiety and mathematics achievement: Evidence from a US sample. Issues in Educational Research, 25(3), 246-259.

Devine, A., Hill, F., Carey, E., & Szűcs, D. (2018). Cognitive and emotional math problems largely dissociate: Prevalence of developmental dyscalculia and mathematics anxiety. Journal of Educational Psychology, 110(3), 431.

Dowker, A., Sarkar, A., & Looi, C. Y. (2016). Mathematics Anxiety: What Have We Learned in 60 Years? *Frontiers in Psychology*, *7*. https://doi.org/10.3389/fpsyg.2016.00508

Dowker, A. (2005). Individual differences in arithmetic: Implications for psychology, neuroscience and education. Psychology Press.

Eysenck, M. W., Derakshan, N., Santos, R., & Calvo, M. G. (2007). Anxiety and cognitive

performance: Attentional control theory. Emotion, 7, 336–353. https://doi.org/10.1037/1528-3542.7.2.336

Foley, A. E., Herts, J. B., Borgonovi, F., Guerriero, S., Levine, S. C., & Beilock, S. L. (2017). The math anxiety-performance link: A global phenomenon. *Current Directions in Psychological Science*, *26*(1), 52–58. https://doi.org/10.1177/0963721416672463

Geary, D. C., Bailey, D. H., Littlefield, A., Wood, P., Hoard, M. K., & Nugent, L. (2009). First-grade predictors of mathematical learning disability: A latent class trajectory analysis. *Cognitive Development*, *24*(4), 411–429. https://doi.org/10.1016/j.cogdev.2009.10.001

Gray, S. A., & Reeve, R. A. (2016). Number-specific and general cognitive markers of preschoolers' math ability profiles. *Journal of Experimental Child Psychology*, *147*, 1–21. https://doi.org/10.1016/j.jecp.2016.02.004

Hart, S. A., & Ganley, C. M. (2019). The nature of math anxiety in adults: Prevalence and correlates. Journal of numerical cognition, 5(2), 122.

Henry, D. L., Baltes, B., & Nistor, N. (2014). Examining the relationship between math scores and English language proficiency. Journal of Educational Research and Practice, 4(1), 2.

Harari, R. R., Vukovic, R. K., & Bailey, S. P. (2013). Mathematics Anxiety in Young Children: An Exploratory Study. *The Journal of Experimental Education*, *81*(4), 538–555. https://doi.org/10.1080/00220973.2012.727888

Hembree, R. (1990). The Nature, Effects, and Relief of Mathematics Anxiety. *Journal for Research in Mathematics Education*, *21*(1), 33–46. JSTOR. https://doi.org/10.2307/749455

Honaker, J., King, G., & Blackwell, M. (2019). *Amelia: A Program for Missing Data* (1.7.6) [Computer software]. https://CRAN.R-project.org/package=Amelia

Jamaludin, A., & Hung, D. (2019). Translational specifications of neural-informed game-based interventions for mathematical cognitive development of low-progress learners: A science of learning approach OER Knowledge Bites, 10, 6-7.

Kellogg, J. S., Hopko, D. R., & Ashcraft, M. H. (1999). The Effects of Time Pressure on Arithmetic Performance. *Journal of Anxiety Disorders*, *13*(6), 591–600. https://doi.org/10.1016/S0887-6185(99)00025-0

Krinzinger, H., Kaufmann, L., & Willmes, K. (2009). Math Anxiety and Math Ability in Early Primary School Years. *Journal of Psychoeducational Assessment*, *27*(3), 206–225. https://doi.org/10.1177/0734282908330583

Lee Swanson, H., Howard, C. B., & Saez, L. (2006). Do different components of working memory underlie different subgroups of reading disabilities?. Journal of learning disabilities, 39(3), 252-269.

Lusardi, A., & Mitchell, O. S. (2011). Financial literacy and retirement planning in the United States (No. 0898–2937). National Bureau of Economic Research.

Langsrud, Ø. (2003). ANOVA for unbalanced data: Use type II instead of type III sums of squares. *Statistics and Computing*, *13*(2), 163–167. https://doi.org/10.1023/A:1023260610025

Lenth, R., Buerkner, P., Herve, M., Love, J., Riebl, H., & Singmann, H. (2020). *emmeans: Estimated Marginal Means, aka Least-Squares Means* (1.5.0) [Computer software]. https://CRAN.R-project.org/package=emmeans

Little, R. J. A. (1988). A Test of Missing Completely at Random for Multivariate Data with Missing Values. *Journal of the American Statistical Association*, *83*(404), 1198–1202. JSTOR. https://doi.org/10.2307/2290157

Mammarella, I. C., Hill, F., Devine, A., Caviola, S., & Szűcs, D. (2015). Math anxiety and developmental dyscalculia: A study on working memory processes. Journal of clinical and experimental neuropsychology, 37(8), 878-887.

Menon, V. (2016). Working memory in children's math learning and its disruption in dyscalculia. Current Opinion in Behavioral Sciences, 10, 125-132.

Menon V, Uddin LQ. Saliency, switching, attention and control: a network model of insula function. Brain Struct Funct. 2010; 214(5–6), 655–667. https://doi.org/10.1007/s00429-010-0262-0 PMID: 20512370

Mutlu, Y. (2019). Math Anxiety in Students with and without Math Learning Difficulties. International Electronic Journal of Elementary Education, 11(5), 471-475.

Mattarella-Micke, A., Mateo, J., Kozak, M. N., Foster, K., & Beilock, S. L. (2011). Choke or thrive? The relation between salivary cortisol and math performance depends on individual differences in working memory and math-anxiety. *Emotion*, *11*(4), 1000–1005. https://doi.org/10.1037/a0023224

McClelland, M. M., Cameron, C. E., Connor, C. M., Farris, C. L., Jewkes, A. M., & Morrison, F. J. (2007). Links between behavioral regulation and preschoolers' literacy, vocabulary, and math skills. *Developmental Psychology*, *43*(4), 947–959. https://doi.org/10.1037/0012-1649.43.4.947

Namkung, J. M., Peng, P., & Lin, X. (2019). The Relation Between Mathematics Anxiety and Mathematics Performance Among School-Aged Students: A Meta-Analysis. *Review of Educational Research*, *89*(3), 459–496. https://doi.org/10.3102/0034654319843494

Ng, J. D., Liew, E., Menon, L., Ling, L. P., & Wang, A. (2007). Learning support for mathematics (lsm): using the 4-pronged intervention approach. Proceedings of the Redesigning Pedagogy: Culture, Knowledge and Understanding Conference [Symposium]. Singapore.

NIE Office of Education Research (2022). iMAGINE: Math Game-based Interventions in Neural-informed education. Retrieved 15 Jan 2022 from iMAGINE: Math Game-based Interventions in Neural-informed education | SingTeach | Education Research for Teachers | Research within Reach (nie.edu.sg)

Organisation for Economic Co-operation and Development. (2013). PISA 2012 results: Excellence through equity: Giving every student the chance to succeed (Volume II). Paris, France: OECD Publishing.

Peng, P., Namkung, J., Barnes, M., & Sun, C. (2016). A meta-analysis of mathematics and working memory: Moderating effects of working memory domain, type of mathematics skill, and sample characteristics. Journal of Educational Psychology, 108(4), 455–473. https://doi.org/10.1037/edu0000079

Pletzer B, Kronbichler M, Nuerk HC, Kerschbaum HH. Mathematics anxiety reduces default mode network deactivation in response to numerical tasks. Front Hum Neurosci. 2015; 9:202. https://doi.org/10.3389/fnhum.2015.00202 PMID: 25954179

Psychological Corporation. (2009). WIAT III: Wechsler Individual Achievement Test. Psychological Corp.

R Core Team. (2020). R: A language and environment for statistical computing (4.0.2) [R Foundation for Statistical Computing]. https://www.r-project.org/

Ramirez, G., Chang, H., Maloney, E. A., Levine, S. C., & Beilock, S. L. (2016). On the relationship between math anxiety and math achievement in early elementary school: The role of problem solving strategies. Journal of Experimental Child Psychology, 141, 83–100. https://doi.org/10.1016/j.jecp.2015.07.014

Ramirez, G., Gunderson, E. A., Levine, S. C., & Beilock, S. L. (2013). Math Anxiety, Working Memory, and Math Achievement in Early Elementary School. Journal of Cognition and Development, 14(2), 187–202. https://doi.org/10.1080/15248372.2012.664593

Ramirez, G., Shaw, S. T., & Maloney, E. A. (2018). Math Anxiety: Past Research, Promising Interventions, and a New Interpretation Framework. Educational Psychologist, 53(3), 145–164. https://doi.org/10.1080/00461520.2018.1447384

Reeve, R., Reynolds, F., Humberstone, J., & Butterworth, B. (2012). Stability and change in markers of core numerical competencies. Journal of Experimental Psychology: General, 141(4), 649–666. https://doi.org/10.1037/a0027520

Richardson, F. C., & Suinn, R. M. (1972). The Mathematics Anxiety Rating Scale: Psychometric data. Journal of Counseling Psychology, 19(6), 551–554. https://doi.org/10.1037/h0033456

Raghubar, K. P., Barnes, M. A., & Hecht, S. A. (2010). Working memory and mathematics: A review of developmental, individual difference, and cognitive approaches. Learning and individual differences, 20(2), 110-122.

Ramirez, G., Gunderson, E. A., Levine, S. C., & Beilock, S. L. (2013). Math anxiety, working memory, and math achievement in early elementary school. Journal of Cognition and Development, 14(2), 187-202.

Rotzer, S., Loenneker, T., Kucian, K., Martin, E., Klaver, P., & Von Aster, M. (2009). Dysfunctional neural network of spatial working memory contributes to developmental dyscalculia. Neuropsychologia, 47(13), 2859-2865.

Rubinsten, O., & Henik, A. (2009). Developmental dyscalculia: Heterogeneity might not mean different mechanisms. Trends in cognitive sciences, 13(2), 92-99.

Schaeffer, M. W., Rozek, C. S., Berkowitz, T., Levine, S. C., & Beilock, S. L. (2018). Disassociating the relation between parents' math anxiety and children's math achievement: Long-term effects of a math app intervention. Journal of Experimental Psychology: General, 147(12), 1782–1790. https://doi.org/10.1037/xge0000490

Szczygieł, M. (2019). How to measure math anxiety in young children? Psychometric properties of the modified Abbreviated Math Anxiety Scale for Elementary Children (mAMAS-E). Polish Psychological Bulletin, 50(4), 303-315. https://doi.org/10.24425/ppb.2019.131003

Szczygiel, M. (2020). Gender, general anxiety, math anxiety and math achievement in early school-age children. Issues in Educational Research, 30(3), 1126-1142.

Smedt, B. D., Taylor, J., Archibald, L., & Ansari, D. (2010). How is phonological processing related to individual differences in children's arithmetic skills? Developmental Science, 13(3), 508–520. https://doi.org/10.1111/j.1467-7687.2009.00897.x

Sorvo, R., Koponen, T., Viholainen, H., Aro, T., Räikkönen, E., Peura, P., Tolvanen, A., & Aro, M. (2019). Development of math anxiety and its longitudinal relationships with arithmetic achievement among primary school children. Learning and Individual Differences, 69, 173–181. https://doi.org/10.1016/j.lindif.2018.12.005

Trezise, K., & Reeve, R. A. (2014). Working memory, worry, and algebraic ability. Journal of Experimental Child Psychology, 121, 120–136. https://doi.org/10.1016/j.jecp.2013.12.001

Van Mier, H. I., Schleepen, T. M., & Van den Berg, F. C. (2019). Gender differences regarding the impact of math anxiety on arithmetic performance in second and fourth graders. Frontiers in psychology, 2690.

Vetter, P., Butterworth, B., & Bahrami, B. (2011). A Candidate for the Attentional Bottleneck: Set-size Specific Modulation of the Right TPJ during Attentive Enumeration. Journal of Cognitive Neuroscience, 23(3), 728–736. https://doi.org/10.1162/jocn.2010.21472

Vukovic, R. K., Kieffer, M. J., Bailey, S. P., & Harari, R. R. (2013). Mathematics anxiety in young children: Concurrent and longitudinal associations with mathematical performance. Contemporary Educational Psychology, 38(1), 1–10. https://doi.org/10.1016/j.cedpsych.2012.09.001

Wechsler, D., Pearson Education, Inc., & Psychological Corporation. (2014). WISC-V: Wechsler Intelligence Scale for Children. NCS Pearson, Inc. : PsychCorp.

Wu, S. S., Barth, M., Amin, H., Malcarne, V., & Menon, V. (2012). Math Anxiety in Second and Third Graders and Its Relation to Mathematics Achievement. Frontiers in Psychology, 3. https://doi.org/10.3389/fpsyg.2012.00162

Yeo, W. T., Tan, A. L., & Jamaludin, A. (2022). 'I feel like I'm fighting fire': Teaching the young and educationally disadvantaged. Teaching and Teacher Education, 113, 103665.

Zhang, J., Zhao, N., & Kong, Q. P. (2019). The Relationship Between Math Anxiety and Math Performance: A Meta-Analytic Investigation. Frontiers in Psychology, 10. https://doi.org/10.3389/fpsyg.2019.01613